\documentclass{amsart}

\usepackage{amsmath,amssymb,amsfonts,graphicx}
\usepackage[all]{xy}
\usepackage{epstopdf}

\newtheorem{theorem}{Theorem}
\newcommand{\bt}{\begin{theorem}}
\newcommand{\et}{\end{theorem}}
\newtheorem{lemma}{Lemma}
\newcommand{\bl}{\begin{lemma}}
\newcommand{\el}{\end{lemma}}
\newtheorem{corollary}{Corollary}
\newcommand{\bc}{\begin{corollary}}
\newcommand{\ec}{\end{corollary}}
\newcommand{\beq}{\begin{equation}}
\newcommand{\eeq}{\end{equation}}
\newcommand{\benum}{\begin{enumerate}}
\newcommand{\eenum}{\end{enumerate}}
\newcommand{\N}{\ensuremath{ \mathbf N }}
\newcommand{\Z}{\ensuremath{\mathbf Z}}
\newcommand{\Q}{\ensuremath{\mathbf Q}}
\newcommand{\R}{\ensuremath{\mathbf R}}
\newcommand{\C}{\ensuremath{\mathbf C}}

\newcommand{\mcd}{\ensuremath{ \mathcal D}}
\newcommand{\mcf}{\ensuremath{ \mathcal F}}

\DeclareMathOperator{\orbit}{\text{orbit}}
\newcommand{\bmat}{\left(\begin{matrix}}
\newcommand{\emat}{\end{matrix}\right)}

\DeclareMathOperator{\qqand}{\qquad\text{and}\qquad}

\title{Forests of complex numbers}
\author{Melvyn B. Nathanson}
\address{Department of Mathematics\\
Lehman College (CUNY)\\
Bronx, NY 10468} 
\email{melvyn.nathanson@lehman.cuny.edu}


\date{\today}

\begin{document}

\begin{abstract}
The Calkin-Wilf tree is an infinite  binary tree whose vertices 
are the positive rational numbers.  
Each number occurs in the tree exactly once and in the form $a/b$, 
where are $a$ and $b$ are relatively prime positive integers.  
In this paper, certain subsemigroups of the modular group are used 
to construct similar trees in the set $\mcd_0$ of positive complex numbers. 
Associated to each semigroup is a forest of trees that partitions $\mcd_0$.   
The  fundamental domain  and the set of cusps of the semigroup are defined 
and  computed.  
\end{abstract}

\maketitle

\section{Forests generated by left-right pairs of matrices}

Let $\N$, $\N_0$, \Z, and \Q\ denote, as usual, the sets of positive integers, 
nonnegative integers, integers, and rational numbers, respectively. 
The Calkin-Wilf tree is a rooted infinite binary tree whose vertices 
are the positive rational numbers.  
The root of the tree is 1, and the generation rule is
\[  
\xymatrix{
& \frac{a}{b} \ar[dl]  \ar[dr] & \\
\frac{a}{a+b} &  & \frac{a+b}{b}
}
\]   
Writing $z = a/b$, we can redraw this as follows:
\[  \xymatrix{
& z  \ar[dl]  \ar[dr] & \\
\frac{z }{z +1} &  & z +1 
}
\]  
Every positive rational number occurs exactly once 
as a vertex in the Calkin-Wilf tree, and the geometry of the tree encodes 
beautiful arithmetical relations between rational numbers   
(Bates, Bunder, and Tognetti~\cite{bate-bund-togn10}, 
Bates and Mansour~\cite{bate-mans11}, 
Calkin and Wilf~\cite{calk-wilf00}, 
Chan~\cite{chan11}, 
Dilcher and Stolarsky~\cite{dilc-stol07}, 
Gibbons, Lester, and Bird~\cite{gibb-lest-bird06},
Han, Masuda, Singh, and Thiel~\cite{han14},
Mallows~\cite{mall11},
Mansour and Shattuck~\cite{mans-shat11}, 
Nathanson~\cite{nath15b,nath15c,nath15d},
Reznick~\cite{rezn90}).

Let 
\[
SL_2(\N_0) = \left\{ \bmat a & b \\ c & d \emat : a,b,c,d \in \N_0 \text{ and } ad-bc=1 \right\}
\]
be the semigroup of $2 \times 2$ matrices with nonnegative integral coordinates 
and determinant 1.  
Every matrix $T =  \bmat a & b \\ c & d \emat \in SL_2(\N_0)$ 
defines a function $z \mapsto T(z)$ by  
\[
T(z) = \frac{az+b}{cz+d}.
\]
For example, if 
\[
L_1 =  \bmat 1 & 0 \\ 1 & 1 \emat  \qqand
R_1 =  \bmat 1 & 1 \\ 0 & 1 \emat
\]
then
\[
L_1(z) = \frac{z}{z+1}\qqand   R_1(z) = z+1.
\]
Thus, the generation rule for the Calkin-Wilf tree can also be presented in the form  
\[ 
\xymatrix{
& z  \ar[dl]  \ar[dr] & \\
L_1(z) &  & R_1(z) 
}
\]

Denote the real and imaginary parts of a complex number $z$ by 
$\Re(z)$ and $\Im(z)$, respectively.  
Let $K$ be a non-real subfield of the complex numbers, that is, 
$K \subseteq \C$ and $K \not\subseteq \R$. 
In this paper we consider the set of ``positive complex numbers''  
\[
\mcd_0 = \{ z = x+yi \in K: x > 0 \text{ and } y >0 \}.
\]

\bt                \label{Forest:theorem:SemigroupAction}
If $z\in \mcd_0$ and $T \in SL_2(\N_0)$, then $T(z) \in \mcd_0$.  
Moreover, the function  $(T,z) \mapsto T(z)$ from 
$SL_2(N_0) \times \mcd_0$ into $\mcd_0$ 
defines a semigroup action of $SL_2(\N_0)$ on the set $\mcd_0$.  
\et

\begin{proof}
If $T =  \bmat a & b \\ c & d \emat \in SL_2(\N_0)$ and if $z = x+yi \in \mcd_0$, then 
$x>0$, $y>0$, and $ad-bc=1$.  A standard calculation gives  
\begin{align*}
T(z) & = \frac{az+b}{cz+d} 
= \frac{(az+b)(c\overline{z}+d)}{|cz+d|^2} \\
%& = \frac{ac|z|^2 +ad z  + bc\overline{z}+bd}{|cz+d|^2} \\
%& = \frac{ac|z|^2 + (ad + bc)x + bd}{|cz+d|^2} +  \frac{(ad - bc)yi}{|cz+d|^2} \\
& = \frac{ac|z|^2 + (ad + bc)x + bd}{|cz+d|^2} +  \frac{yi}{|cz+d|^2}. 
\end{align*}
Because $ac|z|^2 + (ad + bc)x + bd$, $|cz+d|^2$, and $y$ 
are positive real numbers, it follows that $T(z) \in \mcd_0$.
Because $(T_1T_2)z=T_1(T_2z)$ for all $T_1, T_2 \in SL_2(\N_0)$ 
and $z \in \mcd_0$, the function  $(T,z) \mapsto T(z)$ 
from $SL_2(N_0) \times \mcd_0$ into $\mcd_0$ 
defines a semigroup action of $SL_2(\N_0)$ on the set $\mcd_0$.  
This completes the proof.  
\end{proof}

\bl            \label{Forest:lemma:NoFixedPoint}
If $T \in SL_2(\N_0)$ and $T \neq I$, 
then the fractional linear transformation $z\mapsto T(z)$ 
has no fixed points in $\mcd_0$.
\el

\begin{proof}
Let 
$T = \bmat a & b \\ c & d \emat$, 
where $a,b,c,d \in \N_0$ and $ad -bc = 1$.  
If $z$ is a fixed point of $T$, 
then 
\[
\frac{az+b}{cz+d} = z.
\]
Equivalently,
\[
cz^2 + (d-a)z - b= 0.
\]
If $c=0$, then  $a=d=1$ and so $z = T(z) = z+b$, hence $b=0$ and $T = I$, 
which is absurd.
If $c \neq 0$, then 
\begin{align*}
z & = \frac{a-d \pm \sqrt{(d-a)^2 +4bc}}{2c}  \\
%& = \frac{a-d \pm \sqrt{(a+d)^2 -4(ad-bc)}}{2c}  \\
& = \frac{a-d \pm \sqrt{(a+d)^2 -4}}{2c}.
\end{align*}
If $z \in \mcd_0$, then $\Im(z) > 0$ and so $(a+d)^2 - 4 < 0$.  
Equivalently, $a+d = |a+d|<2$.  Because $a, b,c$, and $d$ are nonnegative integers, 
it follows that  $0=ad=1+bc \geq 1$, which is also absurd.  
Therefore, $T$ has no fixed points in $\mcd_0$.  
\end{proof}

Let $L$ and $R$ be non-identity matrices in $SL_2(\N_0)$.
The ordered pair $(L,R)$ will be called a \emph{left-right pair} if 
\[
L(\mcd_0)\cap R(\mcd_0) = \emptyset.
\]
Equivalently, $(L,R)$ is a \emph{left-right pair} if 
\[
L(z_1) \neq R(z_2) \qquad \text{for all $z_1, z_2 \in \mcd_0$.}
\]

\bl      \label{Forest:lemma:LR-free}
If $(L,R)$ is a left-right pair, then the subsemigroup $\langle L,R \rangle$ 
of $SL_2(\N_0)$ generated by $\{L,R\}$ is free.  
\el

\begin{proof} 
The semigroup  $\langle L,R\rangle$  is free if and only if the unique solution 
of the matrix equation
\beq            \label{Forest:LRfree}
T_1 T_2 \cdots T_k = T'_1 T'_2 \cdots T'_{\ell}
\eeq
with $k,\ell \in \N_0$ and $T_i, T'_j \in \{L,R\}$ for all $i \in \{1,2,\ldots, k\}$ and $j \in \{1,2,\ldots, \ell\}$ 
is the trivial solution $k=\ell$ and $T_i = T'_i$ for all  $i \in \{1,2,\ldots, k\}$.  
If there is a nontrivial solution, then there is a minimal solution, 
that is, a matrix identity~\eqref{Forest:LRfree} with $k \geq \ell$ and $k$ minimal.

If $\ell \geq 1$, then for every $z \in \mcd_0$ we have
\[
T_1 T_2 \cdots T_k(z) = T'_1 T'_2 \cdots T'_{\ell}(z)
\]
and so
\[
T_1(z_1) = T'_1(z_2)
\]
where 
\[
z_1 =  T_2 \cdots T_k(z)  \in \mcd_0 \qqand z_2 = T'_2 \cdots T'_{\ell}(z) \in \mcd_0.
\]
Because $(L,R)$ is a left-right pair, it follows that $T_1 = T'_1$, and so 
$T_2 \cdots T_k = T'_2 \cdots T'_{\ell}$, which contradicts the minimality of $k$.  

If $\ell = 0$, then $k \geq 1$ and $T_1 T_2 \cdots T_k = I$.
Let $T_1 = \bmat a & b \\ c & d \emat \in SL_2(\N_0)$.
Then
\[
\bmat d & -b \\ -c & a \emat = T_1^{-1}  = T_2 \cdots T_k \in SL_2(\N_0)
\]
and so $b=c=0$ and $a=d=1$, that is, $T_1 = I \in \{ L, R \}$, which is absurd.
This completes the proof.  
\end{proof}

\bl              \label{Forest:lemma:LvRu-pair}
For positive integers $u$ and $v$, let 
\begin{equation}   \label{Forest:parent-LuRv}
L_u =  \bmat 1 & 0 \\ u & 1 \emat  \qqand
R_v =  \bmat 1 & v \\ 0 & 1 \emat.
\end{equation} 
Then $(L_u,R_v)$ is a left-right pair.
\el

\begin{proof}
If $z = x+yi \in \mcd_0$, then 
\[
L_u(z) = \frac{z}{uz+1} = \frac{u(x^2+y^2)+x}{|uz+1|^2} + \frac{yi}{|uz+1|^2}. 
\]
Because  
\[
u(x^2+y^2)+x < u^2(x^2+y^2) + 2ux+1 = (ux+1)^2+ (uy)^2 = |uz+1|^2
\]
it follows that 
\[
\Re(L_u(z)) < 1. 
\]
Similarly, $R_v(z) = (x+v) + yi$ and so
\[
\Re(R_v(z)) = x+v > v \geq 1.
\]
Therefore,  
\[
\Re(L_u(z_1) ) < 1 < \Re( R_v(z_2))
\]
for all $z_1, z_2 \in \mcd_0$, and so $L(\mcd_0)\cap R(\mcd_0) = \emptyset$.
This completes the proof.
\end{proof}

Lemmas~\ref{Forest:lemma:LR-free} and~\ref{Forest:lemma:LvRu-pair}
imply that the set $\{L_u, R_v\}$ freely generates the semigroup $\langle L_u, R_v\rangle$.   
Note that $L_u = L_1^u$ and $R_v = R_1^v$.  
It is a classical result that $\{L_1, R_1\}$ freely generates $SL_2(\N_0)$, 
and this also proves that $\{L_u, R_v\}$ freely generates $\langle L_u, R_v\rangle$ 
(cf. Nathanson~\cite{nath15c}).

Let $(L,R)$ be a pair of matrices in $SL_2(\N_0)$.
We consider the directed graph $\mcf(L,R)$ whose vertex set is $\mcd_0$ 
and whose edge set is 
\[
\{ (z,L(z)):z\in \mcd_0\} \cup \{ (z,R(z)):z\in \mcd_0\}. 
\]
In this graph, every vertex has outdegree 2:
\begin{equation}   \label{Forest:parent-z}
\xymatrix{
& z  \ar[dl]  \ar[dr] & \\
L(z) &  & R(z) 
}
\end{equation}

We call $L(z)$ the \emph{left child} of $z$ and $R(z)$ the \emph{right child} of $z$, 
and we call $z$ the \emph{parent} of $L(z)$ and of $R(z)$.
Because $L$ and $R$ are invertible matrices, if $L(z_1) = L(z_2)$ or if 
$R(z_1) = R(z_2)$ for some $z_1,z_2 \in \mcd_0$, then $z_1=z_2$.
If $(L, R)$ is a right-left pair, then there do not exist $z_1,z_2 \in \mcd_0$ 
such that $L(z_1) = R(z_2) = z$, and so the indegree of $z$ is either 0 or 1.
We call $z$ an \emph{orphan} if it has no parent, that is, if  $z$ has indegree 0.

For example, let $(L_1,R_1)$ be the pair of matrices 
defined by~\eqref{Forest:parent-LuRv} with $u = v =1$.  
The first three generations of descendants of the complex number $z \in \mcd_0$ are:
\[ 
\xymatrix@=0.16cm{
& & & & z  \ar[ddll]  \ar[ddrr] & & & & \\
& & & & & & & & \\
& & \frac{z }{z +1} \ar[ddl]  \ar[ddr]  &  & & & z +1 \ar[ddl]  \ar[ddr]  & & \\
& & & & & & & & \\
&\frac{z }{2z +1}    \ar[ddl]  \ar[ddr]  & & \frac{2z +1}{z +1}  \ar[ddl]  \ar[ddr] &  &  \frac{z +1}{z +2}  \ar[ddl]  \ar[ddr] &  &z +2 \ar[ddl]  \ar[ddr] & \\
& & & & & & & & \\
\frac{z }{3z +1} & &  \frac{3z +1}{2z +1} \quad \frac{2z +1}{3z +2} &  & \frac{3z +2}{z +1} \quad  \frac{z +1}{2z +3} &  & 
\frac{2z +3}{z +2} \quad  \frac{z +2}{z +3} & & z +3 
}
\]  
\\
If $z$ is a variable, then the vertices of this tree are the linear fractional transformations 
associated with the semigroup $SL_2(\N_0)$.  
Nathanson~\cite{nath15b} described some remarkable arithmetical properties of this tree.

\bt
For every left-right pair $(L,R)$ of matrices in $SL_2(\N_0)$, 
the directed graph $\mcf(L,R)$ 
is a forest of infinite binary trees.  
\et

\begin{proof}
We must prove that every connected component of the graph is a tree.  
If not, then some component of the graph contains an undirected cycle 
of length $n \geq 2$, 
that is, a sequence of $n \geq 2$ vertices 
$z_0,z_1,\ldots, z_{n-1},z_n$ such that 
\benum
\item[(i)]
$z_i \neq z_j$ for $0 \leq i < j \leq n-1$ and $z_n = z_0$;
\item[(ii)]
for $i= 0, 1,\ldots, n-1$, either $(z_i,z_{i+1})$ is an edge 
or $(z_{i+1}, z_i)$ is an edge.
\eenum
Suppose that $(z_0,z_1)$ is an edge.  Then $z_0$ is the parent of $z_1$.  
Because $(L,R)$ is a left-right pair, every vertex has at most one parent.  
This implies that $z_1$ is the parent of $z_2$, and so $z_2$ is the parent of $z_3$.
Continuing inductively, we conclude that $z_i$ is the parent of $z_{i+1}$ 
for $i=0,1,\ldots, n-1$.   
Thus, $(z_0,z_1,\ldots, z_n)$ is not only a cycle in the graph, but is a directed cycle.
It follows that there is a sequence of matrices $T_0, T_1,\ldots, T_{n-1}$ 
such that $T_i \in \{L,R\}$ and $T_iz_i = z_{i+1}$ for all $i=0,1,\ldots, n-1$.  
Thus, $T = T_{n-1}\cdots T_1T_0 \in SL_2(\N_0)$ and 
$T(z_0) = z_0$, that is, $z_0 \in \mcd_0$ is a fixed point of $T$.  
Lemma~\ref{Forest:lemma:NoFixedPoint} implies that  $T=I$ 
and so 
\[
T_0^{-1} = T_{n-1}\cdots T_1 \in SL_2(\N_0)
\]
which is impossible because 
(as observed in the proof of Lemma~\ref{Forest:lemma:LR-free})
the only invertible matrix in 
$SL_2(\N_0)$ is the identity matrix $I$.  
The same argument applies if $(z_1, z_0)$ is an edge.  
Thus, every component of the directed graph $\mcf(L,R)$ is a tree, 
and so $\mcf(L,R)$ is a forest.  
Because every vertex in $\mcf(L,R)$ has outdegree 2, it follows that every component 
of $\mcf(L,R)$ is an infinite binary tree.
\end{proof}

Let $(L,R)$ be a left-right pair, and let $\langle L,R\rangle$ be the subsemigroup  
of $SL_2(\N_0)$ generated by $\{L,R\}$.  
Every complex number $w \in \mcd_0$ is the root of an infinite binary tree whose 
vertices are the complex numbers in $\mcd_0$ constructed from the generation rule~\eqref{Forest:parent-z}.  
The \emph{orbit} of $w$, denoted $\orbit(w)$,  is the set of vertices in this tree.
Equivalently, 
\[
\orbit(w) = \{T(w): T \in \langle L, R \rangle \}.
\]

\bl
Let $(L,R)$ be a left-right pair, and consider the forest $\mcf(L,R)$.  
Let $w_1, w_2 \in \mcd_0$.  If $\orbit(w_1) \cap \orbit(w_2) \neq \emptyset$, 
then either $\orbit(w_1) \subseteq \orbit(w_2)$ 
or $\orbit(w_2) \subseteq \orbit(w_1)$.
\el

\begin{proof}
If $\orbit(w_1) \cap \orbit(w_2)  \neq \emptyset$, then there exist matrices 
$T, T' \in \langle L,R\rangle$ such that 
\[
T(w_1) = T'(w_2).
\]
There are sequences of matrices $(T_i)_{i=1}^k$ and $(T'_j)_{j=1}^{\ell}$ 
such that $T_i, T'_j  \in \{L,R \}$ for $i=1,\ldots, k$ and $j = 1,\ldots, \ell$,
with $T = T_1 T_2\cdots T_k$ and $T' = T'_1 T_2\cdots T'_{\ell}$.
Therefore, 
\[
T_1\left(  T_2\cdots T_k(w_1)  \right) = T(w_1) = T'(w_2) 
= T'_1\left(  T'_2\cdots T'_{\ell}(w_2)  \right)
\]
Because $(L,R)$ is a left-right pair, it follows that  $T_1 = T'_1$ 
and so $ T_2\cdots T_k(w_1) =  T'_2\cdots T'_{\ell}(w_2)$.  
Suppose that $k \geq \ell$.  
Continuing inductively, we obtain 
$ T_{\ell+1}\cdots T_k(w_1) =  w_2$, and so $w_2 \in \orbit(w_1)$.
It follows that $\orbit(w_2) \subseteq \orbit(w_1)$.  
This completes the proof.  
\end{proof}

Let $(L,R)$ be a left-right pair, and let $\mcf(L,R)$ be the 
associated forest whose vertices are the  complex numbers in $\mcd_0$.  
Let $w \in \mcd_0$.  
Every element in $\orbit(w) \setminus \{ w\}$ is a \emph{descendant} of $w$, 
and $w$ is an \emph{ancestor} of every element in $\orbit(w) \setminus \{ w\}$.  
An \emph{orphan} is a complex number in $\mcd_0$ with no ancestors.  
A complex number in $\mcd_0$ is a descendent of an orphan if and only if 
it has only finitely many ancestors.  
There is a one-to-one correspondence between the rooted infinite 
binary trees in the forest $\mcf(L,R)$ and the set of orphans.  
The set of orphans, denoted $\Omega(L,R)$, 
is called the \emph{fundamental domain} of the semigroup $\langle L, R \rangle$.  

An \emph{infinite path} in the forest $\mcf(L,R)$ is a sequence 
$(w_n)_{n=1}^{\infty}$ of complex numbers in $\mcd_0$ 
such that $w_{n+1}= L(w_n)$ or $w_{n+1}= R(w_n)$ for all $n \in \N$.
A \emph{cusp} of the semigroup $\langle L, R \rangle$ is the limit of an 
infinite path in  $\mcf(L,R)$.  Thus, $w^*$ is a cusp if 
$w^* = \lim_{n\rightarrow \infty} w_n$, 
where $(w_n)_{n=1}^{\infty}$ is a path in $\mcd_0$.  
Of course, not every infinite path has a limit.  

For every left-right pair  $(L,R)$ of matrices in $SL_2(\N_0)$, 
we have the following problems:
\benum
\item
Compute the fundamental domain $\Omega(L, R)$.  

\item
Determine if the forest $\mcf(L,R)$ contains infinite binary trees without roots, 
and describe them.  

\item
Determine the cusps of the semigroup $\langle L, R \rangle$.  
\eenum

\section{Example: Trees in the forest $\mcf(L_u,R_v)$}

Let $u$ and $v$ be positive integers, and let $(L_u,R_v)$ be the 
left-right pair of matrices defined by~\eqref{Forest:parent-LuRv}.  
The special case $u = v=1$ is a complex version of the Calkin-Wilf tree.

Let $u > 0$.  For every positive integer $n$, we define the open half disk 
\begin{align*}
\mcd_n
&  = \left\{x + yi \in \mcd_0: 
\left( x- \frac{1}{2nu} \right)^2 + y^2 <  \left( \frac{1}{2nu}  \right)^2 \right\} \\
&  = \left\{x + yi \in \mcd_0  : nu(x^2 + y^2) <  x \right\}.
\end{align*}
We have the region 
\[
\mcd_0 \setminus \mcd_1 
= \left\{x + yi \in \mcd_0  :  \frac{x}{u(x^2 + y^2)} \leq 1 \right\} 
\]
and, for $n \geq 1$,  the half crescents 
\[  
\mcd_n \setminus \mcd_{n+1} 
= \left\{x + yi \in \mcd_0  : n < \frac{x}{u(x^2 + y^2)} \leq n+ 1 \right\} \neq \emptyset.
\]
The half disks and half crescents satisfy the relations  
\[  
\mcd_0 \supset \mcd_1 \supset \mcd_2 \supset  \cdots \supset \mcd_n  \supset \mcd_{n+1} \supset \cdots,   
\]
\[   
\bigcap_{n=1}^{\infty} \mcd_n = \emptyset,
\]
and, for  $0 \leq n < m$,  
\[
 (\mcd_n \setminus \mcd_{n+1}) \cap \mcd_m = \emptyset.
\]
If $w$ is in the half disk $\mcd_n$, then 
\[
|w| \leq \left| w-\frac{1}{2nu} \right| + \frac{1}{2nu} < \frac{1}{nu}.
\]

\begin{figure}
   \scalebox{.60}{\includegraphics{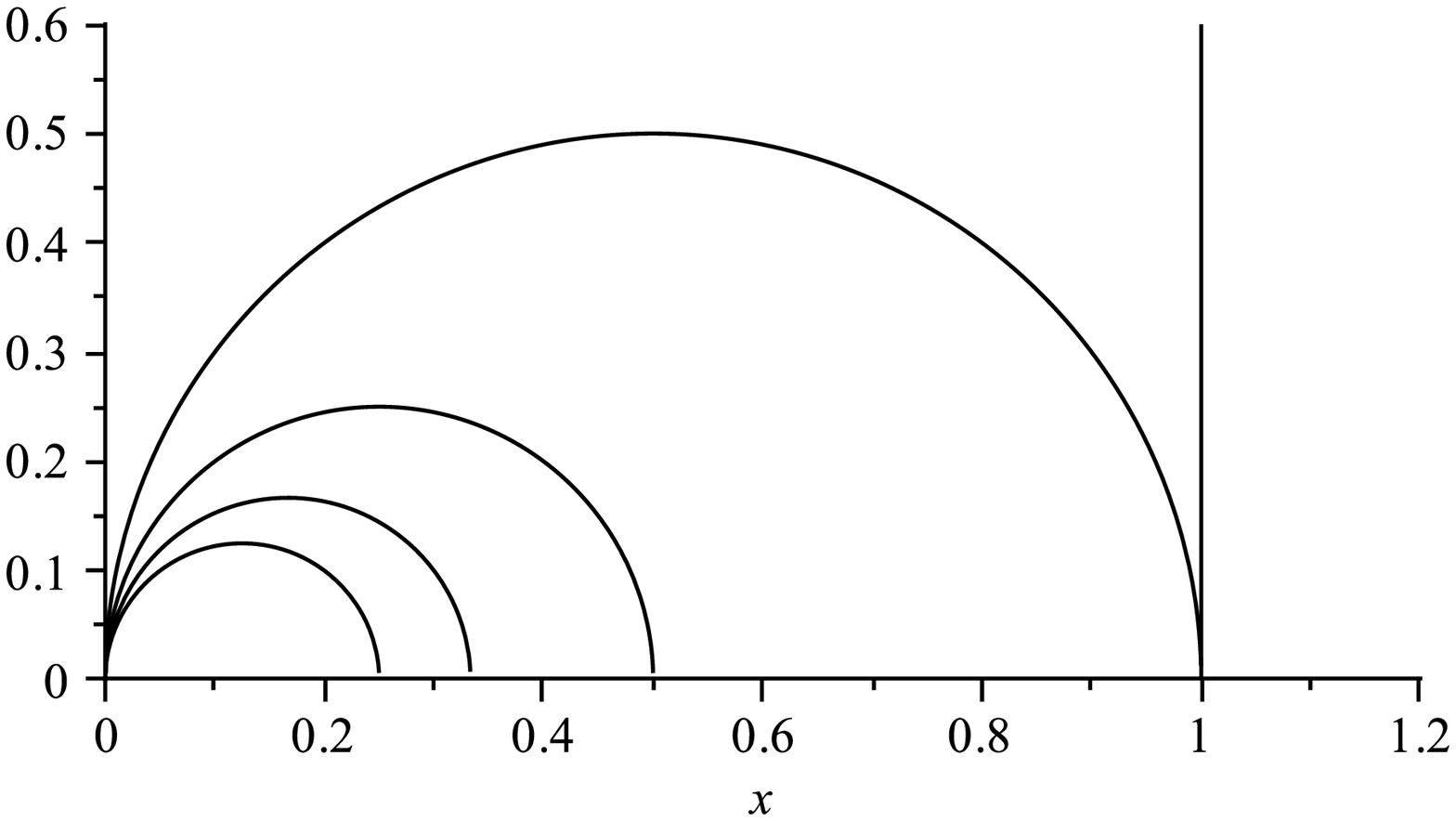}}
\caption{The half crescents for $u=1$ and $n=1,2,3,4$ and the half plane $\Re(z) \geq v=1$.}
\end{figure}

\begin{figure}
   \scalebox{.60}{\includegraphics{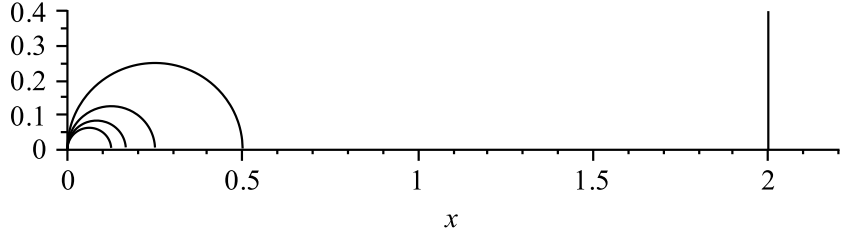}}
\caption{The half crescents for $u=2$ and $n=1,2,3,4$ and the half plane $\Re(z) \geq v=2$.}
\end{figure}

\bl           \label{Forest:lemma:disks}
Let $z\in \mcd_0$, and let $u$, $v$, and $n$ be positive integers.  
Then $\R_v^{-n}(z) \in \mcd_0$ if and only if $\Re(z) > z+nv$, 
and $\L_u^{-n}(z) \in \mcd_0$ if and only if $z \in \mcd_n$.
\el

\begin{proof}
Let $z = x+yi \in \mcd_0$.  
Because $x>0$ and $y > 0$, we have 
\[
R_v^{-n}(z) = z-nv = x-nv + yi  \in \mcd_0
\]
if and only if $\Re(R_v^{-n}(z)) = x-nv > 0$.

Let  $w  = L_u^{-n}(z)$.   
Note that  both $y$ and 
\[
|1-nuz|^2 = (1-nux)^2 + (nuy)^2 
\]
are positive real numbers.  
We have 
\begin{align*} 
w & = L_u^{-n}(z) = \frac{z}{1-nuz} \\ 
& = \frac{z(1-nu\overline{z})}{|1-nuz|^2} = \frac{z - nu|z|^2}{|1-nuz|^2} \\
& = \frac{x - nu(x^2+y^2)}{|1-nuz|^2} + \frac{yi}{|1-nuz|^2} \in \Q(i) 
\end{align*}
and $\Im(w) = y/|1-nuz|^2 > 0$.  
It follows that $w \in \mcd_0$ if and only if  $\Re(w) > 0$ if and only if 
\[
nu(x^2+y^2) < x.
\]
Completing the square, we see that $w \in \mcd_0$ if and only if  
\[
\left( x- \frac{1}{2nu} \right)^2 + y^2 <  \left( \frac{1}{2nu}  \right)^2. 
\]
Thus, $L_u^{-n}(z) \in \mcd_0$ if and only if $z \in D_n$.  
This completes the proof.  
\end{proof}

\bt    \label{Forest:theorem:DiskMaps}
Let $u$ be a positive integer.  The linear fractional transformation  
\[
L_u(z) = \frac{z}{uz+1}
\]
maps  $\mcd_n \setminus \mcd_{n+1}$ onto $\mcd_{n+1} \setminus \mcd_{n+2}$ 
for all  integers $n \geq 0$.   
\et

\begin{proof}
Let $n \geq 0$.
If $z \in \mcd_n \setminus \mcd_{n+1}$ and $w = L_u(z)$, 
then Lemma~\ref{Forest:lemma:disks} implies that 
$L_u^{-n-1}(w) = L_u^{-n}(z)  \in \mcd_0$
and so $w \in \mcd_{n+1}$.  If $w \in \mcd_{n+2}$, then
$L_u^{-n-1}(z)  = L_u^{-n-2}(w) \in \mcd_0$ and so $z \in \mcd_{n+1}$, 
which is absurd.   Therefore, $L_u(z) \in \mcd_{n+1} \setminus \mcd_{n+2}$.

Conversely, let $w \in \mcd_{n+1} \setminus \mcd_{n+2}$, 
and let  $z = L_u^{-1}(w)$.  
Lemma~\ref{Forest:lemma:disks} implies that $L_u^{-n}(z) = L_u^{-n-1}(w) \in \mcd_0$ 
and so $z \in \mcd_n$.  
If $z \in \mcd_{n+1}$, then  $L_u^{-n-2}(w)  = L_u^{-n-1}(z) \in \mcd_0$ and so
$w \in \mcd_{n+2}$, which is absurd.  
Therefore, $z \in \mcd_n \setminus \mcd_{n+1}$, and $L_u(z) = w$.
It follows that the function 
\[
L_u: \mcd_n \setminus \mcd_{n+1} \rightarrow  \mcd_{n+1} \setminus \mcd_{n+2}
\]
is onto.
This completes the proof.  
\end{proof}

\bt
For all  positive integers $u$ and $v$, the fundamental domain 
of the left-right pair $(L_u,R_v)$ is 
\begin{align*}
\Omega(L_u,R_v) & =  \{z \in \mcd_0 : z \notin \mcd_1 \text{ and }  \Re(z) \leq v \} \\
& = \{ x+yi \in \mcd_0 : u(x^2+y^2) \geq x \text{ and } x \leq  v\}.
\end{align*}
\et

\begin{proof}
Let $z=x+yi \in \mcd_0$.   Then $R_v^{-1}(z) = (x-v)+yi \in \mcd_0$ 
if and only if $x > v$.  
Similarly, $L_u^{-1}(z) \in \mcd_0$ if and only if $z \in \mcd_1$.
Thus, the complex number $z$ in $\mcd_0$ has a parent if and only if 
either $\Re(z) > v$ or $z\in \mcd_0$.
Equivalently, $z$ is an orphan if and only if $z \in \mcd_0\setminus \mcd_1$ 
and $\Re(z) > v$.
This completes the proof.  
\end{proof}

Bumby~\cite{bumb14} and Thiel~\cite{thie14} have
independently proved that  every complex number in $\mcd_0$ 
is descended from an orphan with respect to the left-right pair $(L_u,R_v)$.

\bt
For all  positive integers $v$ and $u$, the set of cusps of the semigroup 
$\langle L_u, R_v\rangle$ is $\{0,\infty\}$. 
\et

\begin{proof}
Let $(w_n)_{n=1}^{\infty}$ be an infinite path in the forest $\mcf(L_u,R_v)$.
If $w_{n+1} = R_v(w_n)$  for all $n \geq n_0$, then 
\[
w_n = R_v^{n-n_0}(w_{n_0}) = \Re(w_0)+(n-n_0)v + \Im(w_0)i
\]
and so $\lim_{n\rightarrow \infty} w_n = \infty$.

If $w_{n+1} = L_u(w_n)$  for all $n \geq n_0$, then 
\[
w_n = L_u^{n-n_0}(w_{n_0}) =  \frac{w_{n_0}}{(n-n_0)uw_{n_0}+1}
\] 
and so $\lim_{n\rightarrow \infty} w_n = 0$.

For all $w \in \mcd_0$, we have $L_u(w) \in \mcd_1$ and so
\[
\Re(L_u)(w) < \frac{1}{u} \leq 1.
\]
Similarly, for all $w \in \mcd_0$, we have
\[
R_v(w) - w = v \geq 1.
\]
If $(w_n)_{n=1}^{\infty}$ is an infinite path in the forest $\mcf(L_u,R_v)$ 
such that $w_{n+1} = L_u(w_n)$ for infinitely many $n$ 
and  $w_{n+1} = R_v(w_n)$ for infinitely many $n$, 
then $w_{n+1} \in \mcd_1$ infinitely often
and $w_{n+1} - w_n = v \geq 1$ infinitely often.  
It follows that $\lim_{n\rightarrow \infty} w_n$ does not exist.
Therefore, the set of cusps of the semigroup $\langle L_u, R_v \rangle$ is $\{ 0, 1\}$.
This completes the proof.  
\end{proof}

\section{Open problems}
\benum
\item
Classify the left-right pairs in $SL_2(\N_0)$.

\item
Determine the left-right pairs $(L,R)$ whose associated forests contain 
infinite binary trees without roots.

\item
Let $u, v \in \N$ with $(u,v) \neq (1,1)$.
Find a algorithm  to determine if a matrix belongs to the semigroup   
generated by $L_u$ and $R_v$.

\item
Is there an efficient algorithm  
to determine if two numbers in $\mcd_0$ are in the same tree?

\item
Construct a class of freely generated subsemigroups of $SL_2(\N_))$ of rank $k \geq 3$, 
and describe their associated forests of $k$-regular trees of positive complex numbers.   

\eenum

\def\cprime{$'$} \def\cprime{$'$} \def\cprime{$'$} \def\cprime{$'$}
\providecommand{\bysame}{\leavevmode\hbox to3em{\hrulefill}\thinspace}
\providecommand{\MR}{\relax\ifhmode\unskip\space\fi MR }
\providecommand{\MRhref}[2]{%
  \href{http://www.ams.org/mathscinet-getitem?mr=#1}{#2}
}
\providecommand{\href}[2]{#2}

\end{document}